\documentclass{article}
\usepackage{amsmath,amsthm,amsopn,amstext,amscd,amsfonts,amssymb}
\usepackage{natbib}

\def\tr{\mathop{\rm tr}\nolimits}

\def\etr{\mathop{\rm etr}\nolimits}

\newcommand {\findemo}{\hfill \, $\Box$ \,\\[2ex]}

\renewenvironment{abstract}
                 {\vspace{6pt}
                  \begin{center}
                  \begin{minipage}{5in}
                  \centerline{\textbf{Abstract}}
                  \noindent\ignorespaces
                 }
                 {\end{minipage}\end{center}}
\newtheorem{thm}{\textbf{Theorem}}[section]

\newtheorem{lem}{\textbf{Lemma}}[section]
\theoremstyle{definition}

\title{\Large \textbf{Matrix Kummer-Pearson VII Relation and its Application in Affine Shape}}
\author{
  \textbf{Francisco J. Caro-Lopera} \\
  {\normalsize Department of Basic Sciences} \\
  {\normalsize Universidad de Medell\'{\i}n} \\
  {\normalsize Carrera 87 No.30-65, of. 5-103}\\
  {\normalsize Medell\'{\i}n, Colombia}\\
  {\normalsize E-mail: fjcaro@udem.edu.co}\\[2ex]
  \textbf{Jos\'e A. D\'{\i}az-Garc\'{\i}a} \thanks{Corresponding author\newline
   {\bf Key words.}  Matrix Kummer relation, Pearson VII configuration density, zonal polynomials.\newline
    2000 Mathematical Subject Classification. Primary 62E15; 60E05; secondary
     62H99}\\
  {\normalsize Department of Statistics and Computation} \\
  {\normalsize Universidad Aut\'onoma Agraria Antonio Narro}\\
  {\normalsize 25350 Buenavista, Saltillo, Coahuila, Mexico} \\
  {\normalsize E-mail: jadiaz@uaaan.mx} \\
}
\date{}
\begin{document}
\maketitle
\begin{abstract}
A case of the matrix Kummer relation of \citet{Herz55} based on the Pearson VII type matrix
model is derived in this paper. As a consequence, the polynomial Pearson
VII configuration density is obtained and this set the corresponding exact
inference as a solvable aspect in shape theory.
\end{abstract}

\section{\large Introduction}

Start assuming that a given sample of $n$ ``figures", comprised in $N$ landmarks ("anatomical"
points) in $K$ dimensions, and summarized in an $N\times K$ matrix $\mathbf{X}$, belongs to
certain matrix variate distribution with unknown scale and location parameters. The statistical
theory of shape pursues the distribution of the transforming $\mathbf{X}$ after filtering out
some non important geometrical aspect of the original figures, such as the scale, the position,
the rotation, the uniform shear, and so on; then the so called shape of the population object
should be inferred and some comparisons among population shapes could be performed, among many
others statistical comparisons in the shape space, instead of comparisons in the very noised
non transformed Euclidean space. As usual, the classical theory of shape, with different
geometrical filters, as the Euclidean or affine, was based on Gaussian samples, see for example
\citet{GM93} and the references therein. Then, the non normal applications demanded general
assumptions for the samples, and the so termed generalized shape theory was set under
elliptical models. Under the Euclidean filter we can mention the work of \citet{DC10},
\citet{DC12}; they find the density of all geometrical information about the elliptical random
$\mathbf{X}$ which remains after removing the scale, the position and the rotation. Finally, if
an affine filter is applied, in order to remove from $\mathbf{X}$ all geometrical information
of scale, position and uniform shear, \citet{CDG10} obtained the so term configuration density
of $\mathbf{X}$. All the above densities are expanded in terms of the well known zonal
polynomials, in a series of papers by A.T. James in 60's, see for example \citet{MR1982}.

The transition of the Gaussian shape theory to the elliptical shape theory demanded some
advances in  integration involving zonal polynomials (see for example \citet{CDG10}), but
important problems remain, the computability of the shape densities.

In this paper we focus in alternatives for  such problems. It is easy to see in \citet{GM93},
\citet{DC10} and \citet{DC12} that the structure of shape densities under Euclidean
transformations involves series of zonal polynomials which heritages the difficulties for
computations of the classical hypergeometric series studied by \citet{KE06}. However, a class
of the generalized confluent type series of the configuration densities of \citet{CDG10} can be
handled in order to transformed the series into polynomials, and then the addressed open
problems for computations of the series can be avoided.

The configuration density under Kotz type  samples (including Gaussian) has this property, and
the inference can be performed with polynomials instead of infinite series, see \citet{CDG09}.
The source for this property resides in a generalization of the Kummer relation of
\citet{Herz55}.

This motivates the present work, claiming that there is a similar Kummer relation based on a
Pearson VII distribution, which under certain restriction of the parameters in the associated
Pearson VII configuration density, it can be turned in a polynomial density. Then the inference
can be performed easily by working with the exact likelihood which is written in terms of very
low zonal polynomials. In fact the exact densities can be write down by using formulae for
those polynomials; for example, in the planar shape theory (the most classical applications
resides in the study of figures in $\Re^{2}$), we can use directly the formulae given by
\citet{CDG07} instead of the numerical approaches by \citet{KE06} in order to perform some
analytical properties of the exact density.

This discussion is placed in this paper as follows: section \ref{Sec:KummerPearsonVIIrelation}
defines a Pearson VII type series and finds an integral representation which lead to a matrix
Kummer type relation which we call Kummer-Pearson VII relation; then by applying some general
properties studied by \citet{Herz55}, the equality is extended for the required domains in the
shape theory context. Finally, Kummer-Pearson VII relation gives the finiteness of the Pearson
VII configuration density in section \ref{sec:finitePearsonCD}.

\section{\large Matrix Kummer-Pearson VII relation}

\label{Sec:KummerPearsonVIIrelation} Recall  that the matrix Kummer
relation (due to \citet{Herz55}) states that
\begin{eqnarray}
    {}_{1}F_{1}(a;c;\mathbf{X})&=&\etr(\mathbf{X})\,_{1}F_{1}(c-a;c;-\mathbf{X})
    \quad\mbox{and}\label{eq:Kummer};
\end{eqnarray}
see also \citet{MR1982}.

Now, let $\mathbf{X}>0$ be an $m\times m$ positive definite matrix, then define
\begin{equation}\label{eq:1P1}
  _{1}P_{1}(f(t,\mathbf{X}):a;c;\mathbf{X})=\sum_{t=0}^{\infty}\frac{f(t,\mathbf{X})}{t!}
  \sum_{\tau}\frac{(a)_{\tau}}{(c)_{\tau}}C_{\tau}(\mathbf{X}),
  \end{equation}
where the function $f(t,\mathbf{X})$ is independent of $\tau$, $\tau=(t_{1},\cdots,t_{m})$,
$t_{1}\geq t_{2}\cdots\geq t_{m}>0$, is a partition of $t$,
$$
  (\beta)_{\tau}=\prod_{i=1}^{m}\left(\beta-\frac{1}{2}\left(i-1\right)\right)_{t_{i}}
$$
and
$$
  (b)_{t}=b(b+1)\cdots(b+t-1),\quad (b)_{0}=1.
$$
Then, using this notation we see that the Kummer relation (\ref{eq:Kummer}) is a particular
case of a general type of expressions with the following form
\begin{equation}\label{eq:generalrelation}
    {}_{1}P_{1}\left(f(t,\mathbf{X}):a;c;\mathbf{X}\right)=v(\mathbf{\mathbf{X}})
    {}_{1}P_{1}\left(g(t,\mathbf{X}):c-a;c;h(\mathbf{\mathbf{X}})\right),
\end{equation}
where the functions $v,g$ and $h$ are uniquely determined by the particular function $f$ and
according to the domain of the parameters $a$ and $c$ and the matrix $\mathbf{X}$.

First, we consider an integral representation of the left hand side of
(\ref{eq:generalrelation}) under the model $f(t,\mathbf{X})=(b)_{t}$.
\begin{thm}\label{th:integralrepresentation}
Let $\mathbf{X}<\mathbf{I}$,  $Re(a)>\frac{1}{2}(m-1)$, $Re(c)>\frac{1}{2}(m-1)$ and
$Re(c-a)>\frac{1}{2}(m-1)$. Then for suitable reals  $b$ and $d$, we have that
\begin{eqnarray}\label{eq:integralrepresentation}
    {}_{1}P_{1}\left((b)_{t}d^{-b-t}:a;c;\mathbf{X}\right)=
    \frac{\Gamma_{m}(c)}{\Gamma_{m}(a)\Gamma_{m}(c-a)}\hspace{4.5cm}\nonumber\\
    \times \int_{0<\mathbf{Y}<\mathbf{I}_{m}}\left(d-\tr(\mathbf{XY})\right)^{-b}|\mathbf{Y}|^{a-(m+1)/2}
    |\mathbf{I}-\mathbf{Y}|^{c-a-(m+1)/2}(d\mathbf{Y}).
\end{eqnarray}
\end{thm}
\textit{Proof.} First,  we use a zonal polynomial expansion
\begin{eqnarray*}
    \left(d-\tr(\mathbf{XY})\right)^{-b}&=&\sum_{t=0}^{\infty}\frac{(b)_{t}d^{-b-t}}{t!}[\tr(\mathbf{XY})]^{t}\nonumber\\
    &=&\sum_{t=0}^{\infty}\frac{(b)_{t}d^{-b-t}}{t!}\sum_{\tau}C_{\tau}(\mathbf{XY}).
\end{eqnarray*}
Then integrating term by term using \citet[theorem 7.2.10]{MR1982}, we have that
\begin{eqnarray*}
    &&\int_{0<\mathbf{Y}<\mathbf{I}_{m}}\left(d-\tr(\mathbf{XY})\right)^{-b}|\mathbf{Y}|^{a-(m+1)/2}|\mathbf{I}-
    \mathbf{Y}|^{c-a-(m+1)/2}(d\mathbf{Y})\nonumber\\
    &&\quad=\sum_{t=0}^{\infty}\frac{(b)_{t}d^{-b-t}}{t!}\sum_{\tau}
    \int_{0<\mathbf{Y}<\mathbf{I}_{m}}|\mathbf{Y}|^{a-(m+1)/2}|\mathbf{I}-\mathbf{Y}|^{c-a-(m+1)/2}
    C_{\tau}(\mathbf{XY})(d\mathbf{Y})\nonumber\\
    &&\quad=\sum_{t=0}^{\infty}\frac{(b)_{t}d^{-b-t}}{t!}\sum_{\tau}
    \frac{(a)_{\tau}}{(c)_{\tau}}\frac{\Gamma_{m}(a)\Gamma_{m}(c-a)}{\Gamma_{m}(c)}C_{\tau}(\mathbf{X})\nonumber\\
    &&\quad=\frac{\Gamma_{m}(a)\Gamma_{m}(c-a)}{\Gamma_{m}(c)}
    \sum_{t=0}^{\infty}\frac{(b)_{t}d^{-b-t}}{t!}\sum_{\tau}
    \frac{(a)_{\tau}}{(c)_{\tau}}C_{\tau}(\mathbf{X})\nonumber\\
    &&\quad=\frac{\Gamma_{m}(a)\Gamma_{m}(c-a)}{\Gamma_{m}(c)}\,
    _{1}P_{1}\left((b)_{t}d^{-b-t}:a;c;\mathbf{X}\right),
\end{eqnarray*}
and the required result follows. \findemo Now, we  derive the version of (\ref{eq:Kummer}) but
based on a Pearson VII type model, we call this expression, Kummer-Pearson VII relation.
\begin{thm}\label{th:KummerPearson}
Let $\mathbf{X}>0$,  $Re(a)>\frac{1}{2}(m-1)$, $Re(c)>\frac{1}{2}(m-1)$ and
$Re(c-a)>\frac{1}{2}(m-1)$. Then for suitable reals  $b$ and $d$, the Kummer-Pearson VII
relation is given by
\begin{equation}\label{eq:KummerPearson}
    {}_{1}P_{1}\left((b)_{t}d^{-b-t}:a;c;\mathbf{X}\right)=\left(d-\tr \mathbf{X}\right)^{-b}
    {}_{1}P_{1}\left((b)_{t}\left(d-\tr \mathbf{X}\right)^{-t}:c-a;c;-\mathbf{X}\right).
\end{equation}
\end{thm}
\textit{Proof.} Consider $\mathbf{W}=\mathbf{I}-\mathbf{Y}$ in
(\ref{eq:integralrepresentation}), then we obtain
\begin{eqnarray*}
    &&{}_{1}P_{1}\left((b)_{t}d^{-b-t}:a;c;\mathbf{X}\right)=
    \frac{\Gamma_{m}(c)}{\Gamma_{m}(a)\Gamma_{m}(c-a)}\\&&
    \times \int_{0<\mathbf{W}<\mathbf{I}_{m}}\left(d-\tr[\mathbf{X}(\mathbf{I}-\mathbf{W})]\right)^{-b}
    |\mathbf{W}|^{c-a-(m+1)/2}|\mathbf{I}-\mathbf{W}|^{a-(m+1)/2}(d\mathbf{W})\\&&
    =\frac{\Gamma_{m}(c)}{\Gamma_{m}(a)\Gamma_{m}(c-a)}\\&&
    \times \int_{0<\mathbf{W}<\mathbf{I}_{m}}\left(d-\tr \mathbf{X}-\tr(-\mathbf{XW})\right)^{-b}
    |\mathbf{W}|^{c-a-(m+1)/2}|\mathbf{I}-\mathbf{W}|^{a-(m+1)/2}(d\mathbf{W})\\&&
    =\frac{\Gamma_{m}(c)}{\Gamma_{m}(a)\Gamma_{m}(c-a)}\\&&
    \times \frac{\Gamma_{m}(a)\Gamma_{m}(c-a)}{\Gamma_{m}(c)}\,
    _{1}P_{1}\left((b)_{t}\left(d-\tr \mathbf{X}\right)^{-b-t}:c-a;c;-\mathbf{X}\right),
\end{eqnarray*}
which is the required result.
 \findemo
The reader can compare theorem \ref{th:KummerPearson} (and its proof) with the Kummer relation
(and its proof given by \citet{Herz55}). So the analysis of \citet{Herz55} for extending the
above relations for other values of the parameters, holds in the Kummer-Pearson VII relation
too.

Explicitly, we proved that the integral representation (\ref{eq:integralrepresentation}) of
${}_{1}P_{1}\left((b)_{t}:a;c;\mathbf{X}\right)$ holds for $Re (\mathbf{X})<\mathbf{I}$ (by
analytic continuation), $Re(a)>\frac{1}{2}(m-1)$, $Re(c)>\frac{1}{2}(m-1)$ and
$Re(c-a)>\frac{1}{2}(m-1)$. Then by a suitable modification of the arguments in
\citet{Herz55}, we can extend the domain of  (\ref{eq:KummerPearson}) as follows
\begin{thm}\label{th:KummerPearsonGeneral}
$Re(\mathbf{X})>0$, $Re(a)>\frac{1}{2}(m-1)$ and $Re(c)>\frac{1}{2}(m-1)$. Then for suitable
complex numbers  $b$ and $d$, the Kummer-Pearson VII relation is given by
\begin{equation}\label{eq:KummerPearsonGeneral}
    {}_{1}P_{1}\left((b)_{t}d^{-b-t}:a;c;\mathbf{X}\right)=\left(d-\tr \mathbf{X}\right)^{-b}
    {}_{1}P_{1}\left((b)_{t}\left(d-\tr \mathbf{X}\right)^{-t}:c-a;c;-\mathbf{X}\right).
\end{equation}
\end{thm}

The above relation is important in shape theory applications, in the so called \emph{polynomial
Pearson VII configuration density}.

\section{\large Polynomial Pearson VII Configuration Density}\label{sec:finitePearsonCD}

Our motivation for studying finite shape densities, comes from the computations of
hypergeometric series type involved in these distributions. It is known, that the zonal
polynomials are computable very fast by \citet{KE06}, but the problem now resides in the
convergence and the truncation of the series of zonal polynomials. In fact, in the same
reference of \citet{KE06} we read:

``Several problems remain open, among them automatic detection of convergence .... and it is
unclear how to tell when convergence sets in. Another open problem is to determine the best way
to truncate the series. "

Thus the implicit numerical difficulties for truncation of any configuration density motivate
two areas of investigation: first, continue the numerical approach started by (\citet{KE06}) with
the confluent hypergeometric functions and extend it to the case of some configuration series
type Kotz, Pearson VII, Bessel, Logistic, for example; or second, propose a theoretical
approach for solving   the problem analytically (see \citet{CDG09}).

We study now  the second question corresponding to the
polynomial Pearson VII configuration density.

Recall that a $p\times n$ random matrix $\mathbf{X}$ is said to have a matrix variate symmetric
Pearson type VII distribution with parameters $s,R\in\Re$, $M: p\times n$,
$\mathbf{\Sigma}:p\times p,\,\Phi:n\times n$ with $R>0$, $s>np/2$, $\mathbf{\Sigma}>0$, and
$\Phi>0$ if its probability density function is
$$
\frac{\Gamma(s)}{(\pi
R)^{np/2}\Gamma\left(s-\frac{np}{2}\right)|\mathbf{\Sigma}|^{n/2}|\Phi|^{p/2}}
\left(1+\frac{\tr
(\mathbf{X}-\mathbf{M})'\mathbf{\Sigma}^{-1}(\mathbf{X}-\mathbf{M})\mathbf{\Phi}^{-1}}{R}\right)^{-s}.
$$
When $s=(np+R)/2$, $\mathbf{X}$ is said to have a matrix variate $t$-distribution with $R$
degrees of freedom. And in this case, if $R=1$, then $\mathbf{X}$ is said to have a matrix
variate Cauchy distribution, see\citet{CDG10}.

Then, by \citet{CDG10} we have that (see \citet{GM93} for the gaussian case),
\begin{lem}\label{lem: pearsonviiCDnonisotropic}
Let be
\begin{equation}\label{Aa}
A=\frac{\Gamma_{K}\left(\frac{N-1}{2}\right)}{\pi^{Kq/2}|\mathbf{\Sigma}|^{\frac{K}{2}}
    |U'\mathbf{\Sigma}^{-1}U|^{\frac{N-1}{2}}
    \Gamma_{K}\left(\frac{K}{2}\right)},\quad a=\frac{N-1}{2}
\end{equation}
\begin{equation}\label{X,c}
\mathbf{X}=\frac{1}{R}U'\mathbf{\Sigma}^{-1}\mathbf{\mu}\mathbf{\mu}'
    \mathbf{\Sigma}^{-1}U(U'\mathbf{\Sigma}^{-1}U)^{-1},\quad c=\frac{K}{2},
\end{equation}
\begin{equation}\label{bd}
b=s-\frac{K(N-1)}{2},\quad d=1+
    \frac{\tr\left(\mathbf{\mu}'\mathbf{\Sigma}^{-1}\mathbf{\mu}\right)}{R}.
\end{equation}
If\, $\mathbf{Y}\sim E_{N-1\times K}(\mathbf{\mu}_{N-1\times K},\mathbf{\Sigma}_{N-1\times
N-1}\otimes \mathbf{I}_{K},h)$, for $\mathbf{\Sigma}>0$, then the  non-isotropic noncentral
Pearson type VII configuration density is given by
\begin{equation}\label{eq: pearsonviiCDnonisotropic}
A\,{}_{1}P_{1}\left((b)_{t}d^{-b-t}:a;c;\mathbf{X}\right),
\end{equation}
where ${}_{1}P_{1}(\cdot)$ has been defined in (\ref{eq:1P1}).
\end{lem}
Unfortunately, the above configuration density with general form
$A\,{}_{1}P_{1}(f(t):a;c;\mathbf{X})$ is an infinite series, given that $a=\frac{N-1}{2}$ and
$c=\frac{K}{2}$ are positive (recall that $N$ is the number of landmarks, $K$ is de dimension
and $N-K-1\geq 1$). So a truncation is needed for performing inference when the modified
algorithms of \citet{KE06} are used.

But the above series can be turned into a  polynomial  if we use the following basic principle
of \citet{CDG09}.

\begin{lem}\label{lem:acfinite}
Let $N-K-1\geq1$ as usual, and consider the definition of ${}_{1}P_{1}(\cdot)$ in
(\ref{eq:1P1}). The  infinite configuration density has the general form
$$
  CD_{1}=w(K,N,\mathbf{X})\,\,_{1}P_{1}\left(f(t,\mathbf{X}):\frac{N-1}{2};\frac{K}{2};\mathbf{X}\right),
$$
for suitable functions: $w(\cdot)$, independent of $t$ and $\tau$ but dependent of $K,N$ and
$\mathbf{X}$; and $f(\cdot)$, independent of $\tau$, but dependent of $t$ and possibly of
$\mathbf{X}$  (it depends on the generator elliptical function, compare with the particular
Pearson VII case of lemma \ref{lem: pearsonviiCDnonisotropic}). Then, according to
(\ref{eq:generalrelation}), if the dimension $K$ is even (odd) and the number of landmarks $N$
is odd (even), respectively, then the equivalent configuration density
$$
  CD_{2}=w(K,N,\mathbf{X})v(\mathbf{X})\,\,_{1}P_{1}\left(g(t,\mathbf{X}):-\left(\frac{N-1}{2}-
  \frac{K}{2}\right);\frac{K}{2};h(\mathbf{\mathbf{X}})\right)
$$
is a polynomial of degree $K\left(\frac{N-1}{2}-\frac{K}{2}\right)$ in the latent roots of the
matrix $\mathbf{X}$ (otherwise the series is infinite); where $v, g$ and $h$ are functions
understood in the context of  (\ref{eq:generalrelation}) and depends on the elliptical
generator function.
\end{lem}
Given an elliptical configuration density $CD_{1}$ indexed by the function $f(\cdot)$ and based in the fact that
$a=\frac{N-1}{2}>0$,$c=\frac{K}{2}>0$, the crucial point here consists of finding an integral
representation valid for $c-a=-\frac{N-K-1}{2}<0$, which will lead to  an equivalent elliptical
configuration density $CD_{2}$ indexed by some function $g(\cdot)$. Then the finiteness of $CD_{2}$
follows from $K$ even (odd) and $N$  odd (even), respectively.

In particular, for the Pearson VII generator function, the referred polynomial density is
provided by applying the new Kummer-Pearson VII relation of theorem \ref{th:KummerPearsonGeneral},
via lemma \ref{lem:acfinite}, in lemma \ref{lem: pearsonviiCDnonisotropic}.

\begin{thm}\label{thm:finitePearsonCDisotropic}
Let be $ A, a, b, c, d$ and $\mathbf{X}$ defined  by
(\ref{Aa})-(\ref{bd}).

If\, $\mathbf{Y}\sim E_{N-1\times K}(\mathbf{\mu}_{N-1\times K},\mathbf{\Sigma}_{N-1\times
N-1}\otimes \mathbf{I}_{K},h)$, $\mathbf{\Sigma}>0$, $K$ is even (odd) and $N$ is odd (even),
respectively,  then the polynomial non-isotropic noncentral Pearson type VII configuration
density is given by
\begin{equation}\label{PolynomialPearsonConfigurationDensity}
    \left(d-\tr \mathbf{X}\right)^{-b}
    {}_{1}P_{1}\left((b)_{t}\left(d-\tr \mathbf{X}\right)^{-t}:c-a;c;-\mathbf{X}\right)
\end{equation}

and it is a  polynomial of degree $K\left(\frac{N-1}{2}-\frac{K}{2}\right)$ in the latent roots
of $\mathbf{X}$.
\end{thm}
\textit{Proof.} The proof is trivial, just start with the infinite configuration density
 (\ref{eq: pearsonviiCDnonisotropic}):
$$
A\,{}_{1}P_{1}\left((b)_{t}d^{-b-t}:a;c;\mathbf{X}\right),
$$
where $A, a, b, c, d$ and $\mathbf{X}$ are given by (\ref{Aa})-(\ref{bd}). Then apply
(\ref{eq:KummerPearsonGeneral}) and the result follows. Note that the finiteness follows from
Lemma (\ref{lem:acfinite}) by noting that $c-a=-\frac{N-K-1}{2}$ is a negative integer, when
$K$ is even (odd) and $N$ is odd (even). \findemo

The principle of lemma \ref{lem:acfinite} is based on a known property of the hypergeometric
series easily extended to series of the type ${}_{1}P_{1}(f(t,\mathbf{X}):a,c;\mathbf{X})$, see
(\ref{eq:1P1}), which states that if $a$ is a negative integer or a negative half integer, the
series vanishes in a polynomial.  Then if we have an application following a confluent
distribution type, it is a polynomial, always that the parameter $a$ accepts the addressed
special domain, otherwise the distribution is a series of zonal polynomials and the open
problems for its computability appears. This last case occurs for example in the general
configuration density, which is an infinite confluent series type because $a=\frac{N-1}{2}>0$,
then it is not a trivial to turn that series into polynomials, because the application does not
allow a negative parameter. Then the lemma \ref{lem:acfinite} gives a special subclass by
selecting the number of landmarks (even, odd) given the dimension (odd, even), otherwise the
configuration density remains a series, then the associated Kummer relation must be obtained in
order to transform the numerator parameter $a$ into $c-a$ which is a negative half integer as
required. So, in the context of shape theory under affine transformations, the study of  Kummer
relation type plays an important role.  It is easy to check that the classical Kummer relation
first derive by \citet{Herz55}, and set in the context of zonal polynomials by \citet{CONST63}
is related with a Gaussian kernel, but if we want to obtain some non Gaussian polynomial
densities (under the explained restrictions of $N$ and $K$), then we need to derive new Kummer
relations, it was the case of certain class of Kotz configuration densities (which includes the
classical Gaussian), it need the derivation of the associated Kummer Kotz relation, see
\citet{CDG09}. Thus, in the case of the Pearson polynomial configuration density the
corresponding Kummer relation was the key point for transforming series.

The above discussion opens related problems in some special topics of matrix variate analysis
involving confluent matrix with special domains for the parameters, or general hypergeometric
series type with more or equal Pochhammer symbols in the numerator,  several  examples of this
situations, which demands new developments for Euler relations and similar ones, can be
inferred from some distributions proposed  in \citet[chapters 8--11]{MR1982}.

\section*{\large Acknowledgment}

The first author was supported by Universidad de
Medell\'{i}n-Colombia grant No. 158.

\end{document}